\documentclass{conm-p-l}%
\usepackage{amsmath}
\usepackage{amsfonts}
\usepackage{amssymb}
\usepackage{graphicx}%
\providecommand{\U}[1]{\protect\rule{.1in}{.1in}}

\theoremstyle{definition}

\theoremstyle{remark}

\numberwithin{equation}{section}

\begin{document}
\title[Superiorization Algorithms Applied to Radiation Treatment Planning]{Feasibility-Seeking and Superiorization Algorithms Applied to Inverse Treatment Planning in Radiation Therapy}
\author[R. Davidi]{Ran Davidi}
\address{Department of Radiation Oncology, Stanford University School of Medicine,
Stanford, CA 94305, USA}
\email{rdavidi@stanford.edu}
\author[Y. Censor]{Yair Censor}
\address{Department of Mathematics, University of Haifa, Mount Carmel, Haifa 3190501, Israel}
\email{yair@math.haifa.ac.il}
\author[R. W. Schulte]{Reinhard W. Schulte}
\address{Department of Radiation Medicine, Loma Linda University Medical Center, Loma
Linda, CA 92354, USA}
\email{rschulte@llu.edu}
\author[S. Geneser]{Sarah Geneser}
\address{Department of Radiation Oncology, University of California San Francisco, San
Francisco, CA 94143, USA}
\email{genesers@radonc.ucsf.edu}
\author[L. Xing]{Lei Xing}
\address{Department of Radiation Oncology, Stanford University School of Medicine,
Stanford, CA 94305, USA}
\email{lei@stanford.edu}

\begin{abstract}
We apply the recently proposed superiorization methodology (SM) to the 
inverse planning problem in radiation therapy. The inverse planning problem 
is represented here as a constrained minimization problem of the total variation 
(TV) of the intensity vector over a large system of linear two-sided inequalities. 
The SM can be viewed conceptually as lying between feasibility-seeking for the constraints
and full-fledged constrained minimization of the objective function subject to
these constraints. It is based on the discovery that many feasibility-seeking
algorithms (of the projection methods variety) are perturbation-resilient, and
can be proactively steered toward a feasible solution of the constraints with
a reduced, thus superiorized, but not necessarily minimal, objective function value.

\end{abstract}
\maketitle
















\section*{December 3, 2013. Revised: January 30, 2014}



\section{Introduction}

Computationally demanding numerical minimization techniques are often used 
in optimizing the treatment plan of different types of intensity-modulated radiation therapy 
(IMRT), for example, in volumetric-modulated arc therapy (VMAT). However, some commonly 
employed objective functions and corresponding minimization techniques are not necessarily the
most appropriate for achieving the desired radiation dose distribution
behavior in the patient. This disconnect occurs because minimal solutions to
some current minimization formulations are not guaranteed to provide the desired
dose coverage, conformality, and homogeneity. Therefore, the considerable
computational cost associated with some of these minimization techniques may
not be justified.

We propose to apply the recently developed novel superiorization method (SM)
that improves computational tractability by aiming at a solution that is
guaranteed to satisfy the IMRT planning constraints and results in a
reduced, but not necessarily minimal, value of the objective function. 

The SM can be viewed conceptually as lying between feasibility-seeking for the
constraints and full-fledged constrained minimization of the objective
function subject to these constraints. It is based on the discovery that many
feasibility-seeking algorithms (of the projection methods variety) are
perturbation-resilient, and can be proactively steered toward a feasible
solution of the constraints with a reduced, but not necessarily minimal,
objective function value.

The SM is, thus, capable of producing \textquotedblleft superior feasible
solutions\textquotedblright\ by employing less-demanding feasibility-seeking
projection methods. Therefore, it may replace current computationally demanding 
constrained minimization methods, and potentially
lead to shorter computational times and improved dose distributions.

The paper is laid out as follows. In Section \ref{sec:inv-prob} we briefly
acquaint the reader with the inverse problem of radiation therapy treatment
planning and the corresponding mathematical model. In Section \ref{sec:SM}, a
short review of the SM is given, and in Section \ref{sec:comp}, we present 
an illustrative example that shows how SM can be used to plan a prostate cancer 
IMRT case. Finally, in Section \ref{sec:conc} we provide our conclusions.

\section{The inverse problem of radiation therapy treatment
planning\label{sec:inv-prob}}

Inverse planning is at the heart of intensity-modulated treatment procedures
and critically determines the quality of the resulting treatment plan. Usually, 
the attending radiation oncologist defines the planning target volumes (PTV) 
and the organs at risk (OAR), prescribes the minimum 
and maximum target doses, threshold doses and/or volumes not to be exceeded 
in OAR, and gives importance factors for each. These constraints give rise to 
an inverse problem. 
A solution method is run to find a treatment plan consisting of intensities 
and timing of beam apertures that produce a clinically acceptable dose distribution. 

However, as practiced now, the therapeutic capacity of these applications is
underutilized because of the computing performance of some of the currently
used minimization methods. In this work, we suggest to use the SM to reach an
acceptable treatment plan. Let us first briefly describe the inverse problem
at hand; for more technical details related to different types of IMRT, the reader 
may consult review articles, such as, \cite{Otto2008,Webb2003, Bortfeld2006, Yu1995, CX2002}, 
to name but a few.

IMRT-type techniques are currently the most advanced form of external radiation therapy. 
Similar to its predecessor, 3D conformal radiation therapy (3DCRT), the physician 
must clearly define the objective of the treatment plan by specifying dose and/or 
volume constraints for the PTV and OARs that aim at maximum 
tumor cell killing while minimizing damage to the patient's normal tissues. Whereas 3DCRT 
uses static apertures, the treatment plan resulting from solving the corresponding IMRT 
problem is composed of multiple field directions and the movement of computer controlled 
pairs of multileaf collimator (MLC) leaves for each treatment angle.

The MLC leaves dynamically change during treatment and
modulate the beam to achieve the objectives of the physician-defined 
treatment plan. The beam, therefore, can be conceptually subdivided into a 
two-dimensional grid of beam subunits called beamlets. Finding a clinically acceptable
treatment plan comprised of beam
apertures and weights for the multiple directions and possible locations of
the MLCs is the goal of the inverse treatment planning problem. In the next
paragraph, we discuss a typical model for the inverse treatment planning
problem that leads to a constrained minimization problem, which in turn, fits
the SM framework.

Denote the physician-prescribed dose to the PTV by a dose
vector $d=(d_{j})_{j=1}^{J}\in R^{J}$, where $d_{j}$ is the dose in voxel $j$
of the fully-discretized patient's cross-section. The dose distribution $d$ is
known to have a linear relationship with the intensities of the beamlets,
denoted by an intensity vector $x=(x_{i})_{i=1}^{I}\in R^{I}$, such that
$x_{i}$ is the intensity of the beamlet $i$. The dose computation problem is
formulated as a linear system of equations%
\begin{equation}
Ax=d,\label{eq:sys}%
\end{equation}
where $A$ is the $J\times I$ dose-influence matrix that, when multiplied with the beamlet 
intensity vector, $x$, computes the dose, $d$, at voxels in the patient anatomy. Here, $I$ 
is the total number of beamlets, and $J$ is the total number of voxels.

Further assume that there are $S$ structures (PTV and OAR),
for $s=1,2,\ldots,S,$ and let $O_{s}$ be the set of voxel indices that belong
to each structure, $s$,%
\begin{equation}
O_{s}=\{j_{s,1},j_{s,2,}\ldots j_{s,m(s)}\},
\end{equation}
where $m(s)$ is the number of voxels in the $s$ structure. Then the system
matrix $A$ can be partitioned into blocks%
\begin{equation}
A=\left(
\begin{array}
[c]{c}%
A_{1}\\
A_{2}\\
\vdots\\
A_{S}%
\end{array}
\right)  ,\label{At}%
\end{equation}
so that a submatrix $A_{s}$ will contain the rows of $A$ whose indices appear
in $O_{s},$ and $d_{(s)}$ will be the corresponding subvector of $d$ and the 
system (\ref{eq:sys}) becomes%
\begin{equation}
\left(
\begin{array}
[c]{c}%
A_{1}\\
A_{2}\\
\vdots\\
A_{S}%
\end{array}
\right)  x=\left(
\begin{array}
[c]{c}%
d_{(1)}\\
d_{(2)}\\
\vdots\\
d_{(S)}%
\end{array}
\right)  .
\end{equation}

This typically used method of computing dose does not yet encompass 
the acceptance criteria by which a solution is evaluated by the physician. For treatment planning, the physician is also required to prescribe a target dose for 
each PTV and an upper dose constraint for all OARs. However, the acceptance criteria commonly used to accept or reject a solution are in a dose-volume constraints (DVCs) format. Such criteria specify what percentage part of the structure may deviate from the prescribed dose and by how much (percentage-wise).

Inclusion of such DVCs in the problem model leads to a mixed-integer programming (MIP) optimization problem which for typical clinical case sizes is not easy to solve without resorting to heuristic methods. Attempts to refrain from MIP are not yet well-developed, see, e.g., \cite{Chen2008}.

Following a well-trodden path in this area, with 
roots in \cite{alt84} and \cite{censor-amc}, we replace the system (\ref{eq:sys}) 
by a more flexible model in which the physician specifies lower- and 
upper-dose bounds vectors, $\underline{d}$ and $\overline{d},$ respectively, 
on all voxels in the respective structures. 
For an OAR structure we define:%
\begin{equation}
\overline{d_{(s)}}\equiv d_{(s)},
\end{equation}
and for any target structures $s$ such as the PTV we define:
\begin{equation}
\underline{d_{(s)}}\equiv d_{(s)}.
\end{equation}
Hence, for an OAR we specify:
\begin{equation}
0\leq A_{s}x\leq d_{(s)},\label{eq:oar}%
\end{equation}
and for a target structure, $s$, we require
\begin{equation}
d_{(s)}\leq A_{s}x\leq e_{(s)},\label{eq:target}%
\end{equation}
where $e_{(s)}$ is an additional clinically-specified upper-bound subvector on
the target, which provides a homogeneity constraint for the target dose. Denoting by $a^{t}$ the $t$th row of the matrix $A,$ the
inequalities of (\ref{eq:oar}) are, component-wise,%
\begin{equation}
0\leq\left\langle a^{{\displaystyle j_{s,\ell}}},x\right\rangle \leq
d_{(s)},\text{ for all }\ell=1,2,\ldots,m(s),\text{ }%
\end{equation}
where $j_{s,\ell}\in O_{s}$, for a structure $s,$ and  $\left\langle \cdot,\cdot\right\rangle $ stands for the inner product. The inequalities of
(\ref{eq:target}) are,%
\begin{equation}
d_{(s)}\leq\left\langle a^{{\displaystyle j_{s,\ell}}},x\right\rangle \leq
e_{(s)},\text{ for all }\ell=1,2,\ldots,m(s).
\end{equation}

This leads to a system of linear inequalities%
\begin{equation}
\left(
\begin{array}
[c]{c}%
\underline{d}_{(1)}\\
\underline{d}_{(2)}\\
\vdots\\
\underline{d}_{(S)}%
\end{array}
\right)  \leq\left(
\begin{array}
[c]{c}%
A_{1}\\
A_{2}\\
\vdots\\
A_{S}%
\end{array}
\right)  x\leq\left(
\begin{array}
[c]{c}%
\overline{d}_{(1)}\\
\overline{d}_{(2)}\\
\vdots\\
\overline{d}_{(S)}%
\end{array}
\right)  \label{eq:ineq-system}%
\end{equation}
which serves as the constraints set for the 
minimization problem. For the objective function $\phi$ we use the
\textit{total variation} (TV) of the intensity vector $x,$ given by%
\begin{equation}
\phi\left(  X\right)  =TV(X)=\sum\limits_{u=1}^{U-1}\sum\limits_{u=1}%
^{V-1}\sqrt{\left(  x_{u+1,v}-x_{u,v}\right)  ^{2}+\left(  x_{u,v+1}%
-x_{u,v}\right)  ^{2}},\label{eq:tv}%
\end{equation}
where the two-dimensional array is obtained from the intensity vector $x$ by
$X=\left\{  x_{u,v}\right\}  {}_{u=1,{\text{ }}v=1}^{U,{\text{ }}V}$ where $u$
and $v$ are integers (and $uv=J$). The use of TV minimization in radiation
therapy treatment planning was suggested by Zhu \textit{et al}. in \cite{Zhu2008}, 
which they solved using typical minimization 
approaches, rather than a feasibility problem (\ref{eq:ineq-system}). The TV function
regularizes the objective function
and the inverse problem is formulated as an exact constrained minimization, which results in a
substantial computational burden.

Our approach leads us to the constrained minimization problem
(\ref{eq:mini-prob}) with (\ref{eq:tv}) as the objective and
(\ref{eq:ineq-system}) as the constraints. But instead of attempting to solve this optimization problem we 
use the superiorization methodology described in the next section.

\section{A short review of the SM\label{sec:SM}}

The superiorization methodology (SM) of \cite{jota-compare,med-phys-sm} is
intended for nonlinear constrained minimization (CM) problems of the form%
\begin{equation}
\mathrm{minimize}\left\{  \phi(x)\mid x\in C\right\}  , \label{eq:mini-prob}%
\end{equation}
where $\phi:R^{J}\rightarrow R$ is an objective function and $C\subseteq
\Theta\subseteq{R^{J}}$ is a given feasible set defined by a family of
constraints $\{C_{i}\}_{i=1}^{I},$ where each set $C_{i}$ is a nonempty closed
convex subset of $R^{J},$ so that $C=\cap_{i=1}^{I}C_{i}\neq\emptyset$. Consult 
\cite{jota-compare} and  \cite{med-phys-sm} for details and references on the origins
and development of SM.

In a nutshell, the new paradigm of superiorization lies between
feasibility-seeking and CM. It is not quite trying to solve the full fledged
CM; rather, the task is to find a feasible point that is superior (with
respect to the objective function value) to one returned by a
feasibility-seeking only algorithm.

The SM is beneficial for problems for which an exact CM algorithm has
not yet been discovered, or when existing exact optimization algorithms are
time consuming or require too much computer resources for realistic large problems.
If, in such cases, there exist (space- and time-) efficient iterative
feasibility-seeking projection methods that provide
constraints-compatible solutions, then they can be turned by the SM into
methods that will be practically useful from the point of view of the function
to be optimized. Examples of such situations are given in
\cite{jota-compare,med-phys-sm}.

We associate with the feasible set $C$ a proximity function ${Prox}_{C}%
:\Theta\rightarrow R_{+}$, whose value indicates how incompatible a vector
$x\in\Theta$ is with the constraints. For any given $\varepsilon>0$, a point
$x\in\Theta$ for which ${Prox}_{C}(x)\leq\varepsilon$ is called an
$\varepsilon$\textit{-compatible solution} for $C$. We assume that we have a
feasibility-seeking \textit{algorithmic operator} $\boldsymbol{A}_{C}%
:R^{J}\rightarrow\Theta$, that defines a Basic Algorithm whose iterative step,
given the current iterate vector $x^{k}$, calculates the next iterate
$x^{k+1}$ by%
\begin{equation}
x^{k+1}=\boldsymbol{A}_{C}\left(  x^{k}\right)  . \label{alg:basic}%
\end{equation}
Given $C\subseteq R^{J}$, a proximity function ${Prox}_{C}$, a sequence
$\left\{  x^{k}\right\}  _{k=0}^{\infty}\subset\Theta$ and an $\varepsilon>0,$
then an element $x^{K}$ of the sequence which has the properties: (i)
${Prox}_{C}\left(  x^{K}\right)  \leq\varepsilon,$ and (ii) ${Prox}_{C}\left(
x^{k}\right)  >\varepsilon$ for all $0\leq k<K,$ is called an $\varepsilon
$\texttt{-}\textit{output of the sequence }$\left\{  x^{k}\right\}
_{k=0}^{\infty}$\textit{ with respect to the pair}\texttt{ }$(C,$\texttt{
}${Prox}_{C})$. We denote it by $O\left(  C,\varepsilon,\left\{
x^{k}\right\}  _{k=0}^{\infty}\right)  =x^{K}$, $O$ standing for output.

Clearly, an $\varepsilon$-output $O\left(  C,\varepsilon,\left\{
x^{k}\right\}  _{k=0}^{\infty}\right)  $ of a sequence $\left\{
x^{k}\right\}  _{k=0}^{\infty}$ might or might not exist, but if it does, then
it is unique. If $\left\{  x^{k}\right\}  _{k=0}^{\infty}$ is produced by an
algorithm intended for the feasible set $C,$ such as the Basic Algorithm
(\ref{alg:basic} ), without a termination criterion, then $O\left(
C,\varepsilon,\left\{  x^{k}\right\}  _{k=0}^{\infty}\right)  $ is the
\textit{output} produced by that algorithm when it includes the termination
rule to stop when an $\varepsilon$-compatible solution for $C$ is reached.

In order to \textquotedblleft superiorize\textquotedblright\ such an algorithm
we need it to have \textit{strong perturbation resilience }in the sense that for
every $\varepsilon>0,$ for which an $\varepsilon$-output is defined for a
sequence generated by the Basic Algorithm, for every $x^{0}\in\Theta$, we have
also that the $\varepsilon^{\prime}$-output is defined for every
$\varepsilon^{\prime}>\varepsilon$ and for every sequence $\left\{
y^{k}\right\}  _{k=0}^{\infty}$ generated by $y^{k+1}=\boldsymbol{A}%
_{C}\left(  y^{k}+\beta_{k}v^{k}\right)  ,$ for all $k\geq0,$ where the vector
sequence $\left\{  v^{k}\right\}  _{k=0}^{\infty}$ is bounded and the scalars
$\left\{  \beta_{k}\right\}  _{k=0}^{\infty}$ are such that $\beta_{k}\geq0$,
for all $k\geq0,$ and $\sum_{k=0}^{\infty}\beta_{k}<\infty$. See our recent
\cite{jota-compare} for details.

Along with the constraints set $C\subseteq R^{J}$, we look at an objective
function $\phi:R^{J}\rightarrow R$, with the convention that a point in
$R^{J}$ for which the value of $\phi$ is smaller is considered
\textit{superior} to a point in $R^{J}$ for which the value of $\phi$ is larger.

The essential idea of the SM is to make use of the perturbations in order to
transform a strongly perturbation resilient algorithm that seeks a
constraints-compatible solution for $C$ (i.e., is seeking feasibility) into
one whose outputs are equally good from the point of view of
constraints-compatibility, but are superior (not necessarily optimal)
according to the objective function $\phi$.

This is done by producing from the Basic Algorithm another algorithm, called
its \textit{superiorized} version, that makes sure not only that the
$\beta_{k}v^{k}$ are bounded perturbations, but also that $\phi\left(
y^{k}+\beta_{k}v^{k}\right)  \leq\phi\left(  y^{k}\right)  $, for $k\geq L$
for some integer $L\geq0$. The Superiorized Version of the Basic Algorithm
assumes that we have available a summable sequence $\left\{  \eta_{\ell
}\right\}  _{\ell=0}^{\infty}$ of positive real numbers (for example,
$\eta_{\ell}=a^{\ell}$, where $0<a<1$) and it generates, simultaneously with
the sequence $\left\{  y^{k}\right\}  _{k=0}^{\infty}$ in $\Theta$, sequences
$\left\{  v^{k}\right\}  _{k=0}^{\infty}$ and $\left\{  \beta_{k}\right\}
_{k=0}^{\infty}$. The latter is generated as a subsequence of $\left\{
\eta_{\ell}\right\}  _{\ell=0}^{\infty}$, resulting in a nonnegative summable
sequence $\left\{  \beta_{k}\right\}  _{k=0}^{\infty}$. The algorithm further
depends on a specified initial point $y^{0}\in\Theta$ and on a positive
integer $N$. It makes use of a logical variable called \textit{loop}\emph{.
}The Superiorized Version of the Basic Algorithm is presented next by its pseudo-code.\medskip{}

\textbf{The Superiorized Version of the Basic Algorithm}

\textbf{set} $k=0$

\textbf{set} $y^{k}=y^{0}$

\textbf{set} $\ell=-1$

\textbf{repeat}

$\qquad$\textbf{set} $n=0$

$\qquad$\textbf{set} $y^{k,n}=y^{k}$

$\qquad$\textbf{while }$n$\textbf{$<$}$N$

$\qquad\qquad$\textbf{set }$v^{k,n}$\textbf{ }to be a nonascending vector for
$\phi$ at $y^{k,n}$

$\qquad$\textbf{$\qquad$set} \emph{loop=true}

$\qquad$\textbf{$\qquad$while}\emph{ loop}

$\qquad\qquad\qquad$\textbf{set $\ell=\ell+1$}

$\qquad\qquad\qquad$\textbf{set} $\beta_{k,n}=\eta_{\ell}$

$\qquad\qquad\qquad$\textbf{set} $z=y^{k,n}+\beta_{k,n}v^{k,n}$

$\qquad\qquad\qquad$\textbf{if }$\phi\left(  z\right)  $\textbf{$\leq$}%
$\phi\left(  y^{k}\right)  $\textbf{ then }

$\qquad\qquad\qquad\qquad$\textbf{set }$n$\textbf{$=$}$n+1$

$\qquad\qquad\qquad\qquad$\textbf{set }$y^{k,n}$\textbf{$=$}$z$

$\qquad\qquad\qquad\qquad$\textbf{set }\emph{loop = false}

$\qquad$\textbf{set }$y^{k+1}$\textbf{$=$}$\boldsymbol{A}_{C}\left(
y^{k,N}\right)  $

$\qquad$\textbf{set }$k=k+1$ \medskip{}

$\qquad$Analysis of the Superiorized Version of the Basic Algorithm
\cite{jota-compare,med-phys-sm}, shows that it produces outputs that are
essentially as constraints-compatible as those produced by the original (not
superiorized) Basic Algorithm. However, due to the repeated steering of the
process toward reducing the value of the objective function $\phi$, we can
expect that the output of the Superiorized Version will be superior (from the
point of view of $\phi$) to the output of the original
algorithm.\textbf{\smallskip\ }A recent work that includes results about the
SM appears in this volume \cite{swiss-knife}.

\section{Demonstrative examples\label{sec:comp}}

The anonymized pelvic planning CT (computed tomography) of a prostate
cancer patient was employed for the IMRT treatment planning using our proposed
method. Seven equispaced fields were used for targeting the PTV. The 
acceptance criteria were set using the RTOG 0815 randomized trial protocol
\cite{RTOG0815}. The following dose intervals were chosen empirically and used 
as the lower- and upper- dose constraints: Rectum [0-30], Bladder [0-75], Body [0-5], 
Small Bowel [0-20] and Prostate (PTV) [81-82]. 

Our demonstration of the approach was done by comparing the outputs of a
TV-superiorization algorithm with an, otherwise identical, algorithm that
aimed at only satisfying the dose constraints, without applying the SM. 
Here $\boldsymbol{A}_{C}$ was chosen to be the algebraic reconstruction technique (ART) 
for inequalities \cite{HL1978}. It was proven to be bounded perturbation resilient (although without using this term) in \cite{BDHK}, and strongly perturbation resilient in \cite{med-phys-sm}.

From a radiation delivery stand point, a solution that is easy to deliver is
one that has a piecewise constant intensety-beamlet map. The reason has to do
with the physical constraints coming from the MLCs, they require that the
beamlets have a small number of signal levels. It was, therefore, suggested in
the literature to use total-variation (TV) to force the solution to be
piecewise constant \cite{Block2007, Kolehmainen2006}.

We performed two experiments with different starting conditions. For the first
experiment, we initiated the algorithm with the zero vector of beamlet intensities
and for the second experiment all  beamlet intensities were given the value 10. Tables
1 and 2 summarize the results for the two experiments, and in Figure 1 we
present the associated DVH (dose-volume histogram) curves for the prostate plan.

For the first experiment, the TV-superiorization algorithm produced a solution
that met the acceptance criteria after 12 iterations whereas the
feasibility-seeking algorithm was not able to reach an acceptable solution
after this number of iterations. For the second experiment, the
TV-superiorization algorithm reached an acceptable solution even faster, i.e.,
after 7 iterations, and the feasibility-seeking algorithm again failed some of
the acceptance criteria after this number of iterations.

\begin{table}[th]
\caption{RTOG 0815 acceptance criteria and results of experiment 1 described
in Section 4 (TVS stands for TV-superiorization)}%
\label{eqtable}%
\renewcommand\arraystretch{1.5} \noindent%
\[%
\begin{array}
[c]{|c|c|c|}\hline
\text{Acceptance criteria} & \text{Exp 1 with TVS} & \text{Exp 1 without
TVS}\\\hline
\text{PTV: Min Allowed Dose: 75.24 Gy} & \text{75.24 Gy} & \text{56.13
Gy}\\\hline
\text{PTV: Max Allowed Dose: 84.74 Gy} & \text{84.69 Gy} & \text{89.42
Gy}\\\hline
\text{Rectum: No more than 50\% of the} & \text{34.50 \%} & \text{8.50 \%}\\
\text{volume should exceed 60.00 Gy} & \text{} & \text{}\\\hline
\text{Rectum: Max Dose} & \text{82.64 Gy} & \text{82.71 Gy}\\\hline
\end{array}
\]
\end{table}

\begin{table}[th]
\caption{RTOG 0815 acceptance criteria and results of experiment 2 described
in Section 4 (TVS stands for TV-superiorization)}%
\label{eqtable}%
\renewcommand\arraystretch{1.5} \noindent%
\[%
\begin{array}
[c]{|c|c|c|}\hline
\text{Acceptance criteria} & \text{Exp 2 with TVS} & \text{Exp 2 without
TVS}\\\hline
\text{PTV: Min Allowed Dose: 75.24 Gy} & \text{77.80 Gy} & \text{76.15
Gy}\\\hline
\text{PTV: Max Allowed Dose: 84.74 Gy} & \text{84.71 Gy} & \text{87.63
Gy}\\\hline
\text{Rectum: No more than 50\% of the} & \text{36.90 \%} & \text{40.50 \%}\\
\text{volume should exceed 60.00 Gy} & \text{} & \text{}\\\hline
\text{Rectum: Max Dose} & \text{84.09 Gy} & \text{87.25 Gy}\\\hline
\end{array}
\]
\end{table}


\begin{figure}[tb]
\includegraphics[scale=0.3]{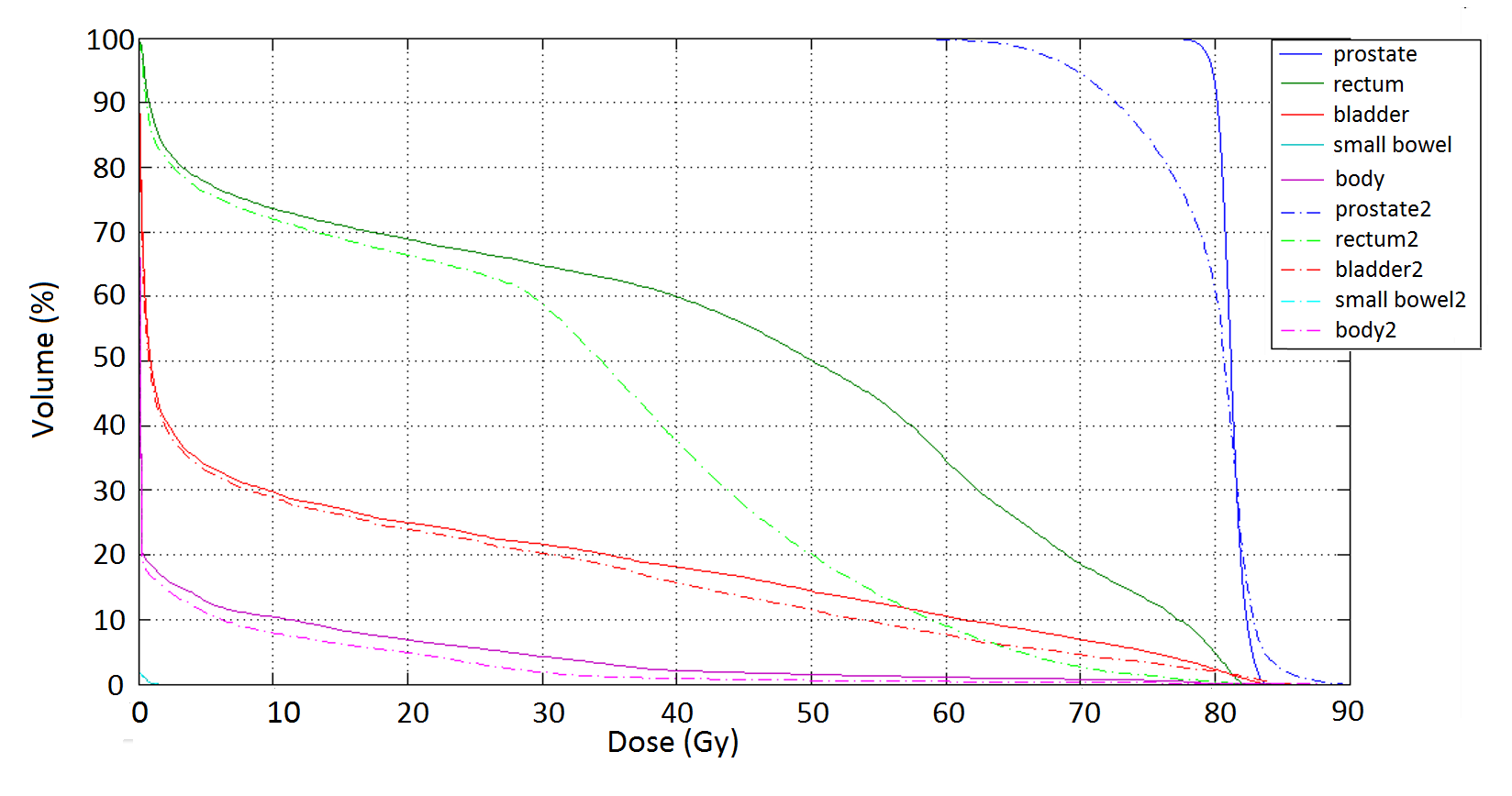}
\includegraphics[scale=0.3]{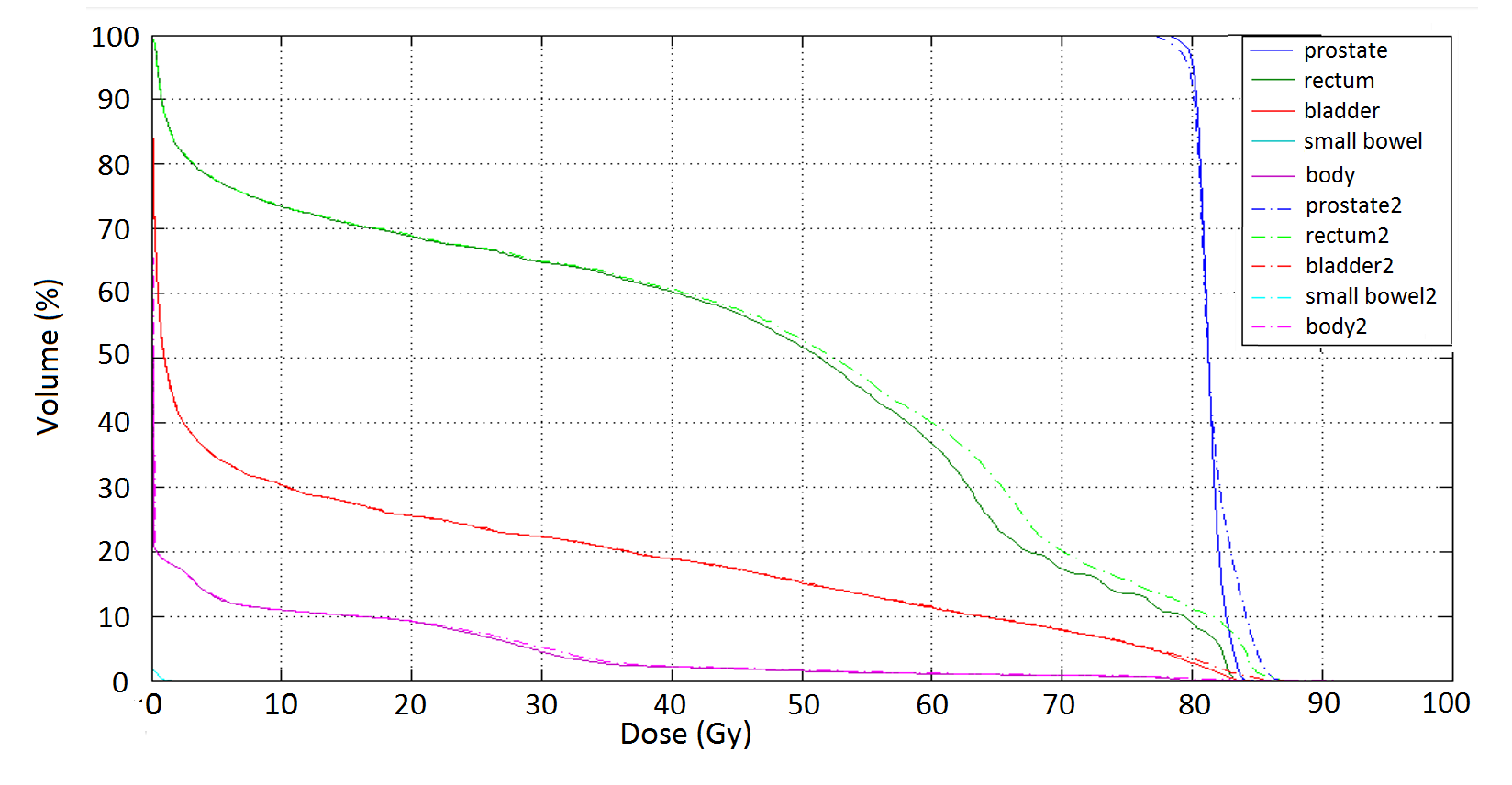}
\caption{Dose volume histograms (DVH) of the two experiments. Solid lines
represent the algorithm with TV-superiorization (broken lines represent the algorithm without TV-superiorization). The first (top) took 12 iterations and the second (bottom)
took 7 iterations.}%
\label{firstfig}%
\end{figure}

\section{Conclusions\label{sec:conc}}
Our proposed method successfully produced conformal solutions that met the acceptance
criteria while an otherwise identical algorithm without superiorization failed to do so with
the same number of iterations. Future work will assess the computational gain of the
superiorization method compared to a conventional method and investigate its utility for
computationally more complex problems that can be found in modulated techniques for arc
therapy.
\bigskip

\textbf{Acknowledgements}. This work is supported by the U.S. Department of
Defense Prostate Cancer Research Program Award No. W81XWH-12-1-0122, by grant
number 2009012 from the United States--Israel Binational Science Foundation
(BSF) and by the U.S. Department of the Army Award No. W81XWH-10-1-0170. Some
of the material was presented as a poster at the Technology for Innovation in
Radiation Oncology, Joint Workshop of the American Society for Radiation
Oncology (ASTRO), the National Cancer Institute (NCI) and the American
Association of Physicists in Medicine (AAPM), National Institutes of Health
(NIH), Bethesda, MD, USA, June 13--14, 2013.

\end{document}